\documentclass[12pt,a4paper]{article}
\usepackage{amsfonts}
\usepackage{amsmath}
\usepackage{makeidx}
\usepackage{graphicx}

\newtheorem{theorem}{Theorem}[section]

\newtheorem{corollary}[theorem]{Corollary}

\newtheorem{definition}{Definition}[section]

\newtheorem{lemma}[theorem]{Lemma}

\newtheorem{problem}{Problem}
\newtheorem{proposition}[theorem]{Proposition}

\newenvironment{proof}[1][Proof]{\textbf{#1.} }{\ \rule{0.5em}{0.5em}}

\begin{document}

\title{Sub-Gaussian short time asymptotics for measure metric Dirichlet
spaces }
\author{Andr\'{a}s Telcs \\
%EndAName
{\small Department of Computer Science and Information Theory, }\\
{\small University of Technology and Economics Budapest}\\
{\small Goldmann Gy\"{o}rgy t\'{e}r 3, V2, 138}\\
{\small Budapest,}\\
{\small H-1111, HUNGARY}\\
{\small telcs@szit.bme.hu}}
\maketitle

\begin{abstract}
This paper presents estimates for the distribution of the exit time from
balls and short time asymptotics for measure metric Dirichlet spaces. The
estimates cover the classical Gaussian case, the sub-diffusive case which
can be observed on particular fractals and further less regular cases as
well. The proof is based on a new chaining argument and it is free of volume
growth assumptions.

MSC2000 31C05, 60J45, 60J60
\end{abstract}

\section{Introduction}

The short-time asymptotics of the heat kernel for Riemannian manifolds has
the classical form due to Varadhan \cite{V}:
\begin{equation*}
\lim_{t\rightarrow 0}t\log p_{t}\left( x,y\right) =-\frac{1}{4}d^{2}\left(
x,y\right) .
\end{equation*}%
It has got recently a lot of attention generalizing the short-time
asymptotics to Dirichlet spaces. Such type of results obtained in Ram\'{\i}%
rez \cite{Ra1} and Hino, Ram\'{\i}rez \cite{HR}( see also Norris \cite{N},
Sturm \cite{S}). The use of the intrinsic metric provides Gaussian
estimates, $d^{2}\left( x,y\right) $ appear as in the case of $%
%TCIMACRO{\U{211d} }%
%BeginExpansion
\mathbb{R}
%EndExpansion
^{d}$ or Riemannian manifolds. \ A vast amount of papers was devoted to
explore further properties of the heat kernel on Riemannian manifolds. \ \
Necessary and sufficient condition where provided for the local two-sided
Gaussian estimate
\begin{equation}
\frac{c}{V(x,\sqrt{t})}\exp \left( -\frac{d^{2}(x,y)}{ct}\right) \leq
p_{t}(x,y)\leq \frac{1}{cV(x,\sqrt{t})}\exp \left( -\frac{cd^{2}(x,y)}{t}%
\right) .  \label{LiYau}
\end{equation}%
of the heat kernel $p_{t}\left( x,y\right) $ \ on Riemannian manifolds in
Grigor'yan \cite{Gr1}, Saloff-Coste \cite{SC}. Here $V\left( x,r\right) $
stands for the volume (with respect. of the measure given on the space) of
the ball $B\left( x,r\right) $ centered at $x$ with radius $r.$ Similar but
sub-Gaussian upper- and two-sided estimates have been obtained for
particular fractals (see Barlow \cite{B2}). \ In \cite{GT2}\ sufficient and
necessary conditions were given for the sub-Gaussian estimates
\begin{equation}
p_{t}(x,y)\geq \frac{1}{CV(x,t^{1/\beta })}\exp \left( -C\left( \frac{%
d^{\beta }(x,y)}{t}\right) ^{\frac{1}{\beta -1}}\right)  \label{LEb}
\end{equation}%
\begin{equation}
p_{t}(x,y)\leq \frac{1}{cV(x,t^{1/\beta })}\exp \left( -c\left( \frac{%
d^{\beta }(x,y)}{t}\right) ^{\frac{1}{\beta -1}}\right)  \label{UEb}
\end{equation}%
for weighted graphs and in \cite{GT3} for measure metric spaces.

\ Consider $T_{x,R}$, the exit time from a ball $B\left( x,R\right) $. If
the stating point $X_{0}=y$ let us denote the expected value of $T_{B\left(
x,R\right) }$ by $\mathbb{E}_{y}\left( x,R\right) $ \ and by $E\left(
x,R\right) $ \ if $y=x.$ Let us denote by $a_{i}\simeq b_{i}$ the fact that
there is a $c>0$ such that $\frac{1}{C}<\frac{a_{i}}{b_{i}}<C$ for all $i$.
On many fractals (or fractal type graph) there is a space-time scaling
function satisfying $\left( E_{F}\right) :$
\begin{equation}
E\left( x,R\right) \simeq F\left( R\right) ,
\end{equation}%
in particular $\left( E_{\beta }\right) :$%
\begin{equation}
E\left( x,R\right) \simeq R^{\beta },  \label{Ebeta}
\end{equation}%
for a $\beta \geq 2.$

During the proof of the upper estimate an interesting side-result can be
observed $\left( E_{\beta }\right) $ implies that
\begin{equation}
\mathbb{P}\left( T_{x,R}<t|X_{0}=x\right) \leq C\exp \left( -\left( \frac{%
R^{\beta }}{Ct}\right) ^{\frac{1}{\beta -1}}\right) .  \label{PUE}
\end{equation}%
\

One might wonder about the conditions which ensure similar lower bound or
small time asymptotics for the heat kernel. The proof of the upper bound $%
\left( \ref{PUE}\right) $ uses some kind of chaining argument (cf. \cite{B2}
and \cite{GT2}). \ The off-diagonal lower estimates are typically shown
using the chaining argument\ by Aronson \cite{A} which uses a volume growth
condition. In the present paper we provide lower counterpart of $\left( \ref%
{PUE}\right) $ and the short time asymptotics based on a new chaining
argument ( Proposition \ref{p3} ). This chaining argument provides a lower
estimate for the distribution of the hitting time of a ball, which we think
might has some interest on its own. No volume growth or bounded covering
conditions are needed.

In order to weaken the restriction on the mean exit time let us introduce
some notions.

\begin{definition}
\label{dKL}The sub-Gaussian kernels defined as follows. Let $A\subset
M,t.R>0 $ and $\kappa =\kappa _{A}\left( t,R\right) $ is the largest integer
for which
\begin{equation*}
\frac{t}{\kappa }\leq q\inf\limits_{y\in A}E\left( y,\frac{R}{\kappa }\right)
\end{equation*}%
in particular,%
\begin{equation*}
\kappa \left( x,t,R\right) =\kappa _{B\left( x,R\right) }\left( t,R\right)
\end{equation*}%
and similarly $\nu =\nu _{A}\left( t,R\right) $ \ is the smallest for which
\begin{equation*}
\frac{t}{\nu }\geq q\sup\limits_{y\in A}E\left( y,\frac{R}{\nu }\right) ,
\end{equation*}%
\begin{equation*}
\nu \left( x,t,R\right) =\nu _{B\left( x,R\right) }\left( x,t,R\right) .
\end{equation*}%
We define $\kappa =0$ and $\nu =\infty $ if there is no such an integer. The
constants $q>0$ will be specified later.
\end{definition}

These integers will define the number of iterations we will use in chaining
arguments which lead to the upper and lower estimates.

Our approach is different from those one \cite{HR}-\cite{N} which use the
intrinsic metric and recapture the usual Gaussian $R\rightarrow R^{2}$
scaling. \ We assume that the metric is given, predefined and we obtain a
picture of the heat diffusion with respect to this metric. \ As Ramirez in
\cite{Ra1} points out the two approaches are complement each other.

In the whole sequel we consider $\left( M,\mu ,d\right) $ a locally compact
separable measure metric space with Radon measure $\mu ,$ with full
support.\ The metric is assumed to be a geodesic one. \ A strictly local
regular Dirichlet form $\left( \mathcal{E},\mathcal{F}\right) $ in $%
L^{2}\left( M,\mu \right) $ is considered and let $\left( X_{t}\right) $ be
the associated diffusion process on $M$ (cf. \cite{FOT}). The corresponding
Feller semigroup is $P_{t}.$ Denote $\mathbb{P}_{x},\mathbb{E}_{x}$ the
probability measure and expectation given for $X_{0}=x\in M.$ We assume that
$\left( X_{t}\right) $ has a transition density \ $p_{t}\left( x,y\right) $
with respect to $\mu ,$ furthermore $p_{t}\left( x,y\right) $ satisfies the
following property:

$\left( 1.\right) $ $\ p_{t}\left( x,y\right) \geq 0,$

$\left( 2.\right) $ $\int_{M}p_{t}\left( x,y\right) d\mu \left( y\right) =1,$

$\left( 3.\right) $ $p_{t}\left( x,y\right) =p_{t}\left( y,x\right) ,$

$\left( 4.\right) $ $p_{t}\left( x,y\right) =\int p_{s}\left( x,z\right)
p_{t-s}\left( z,y\right) d\mu \left( z\right) .$

Now we give the definition of the elliptic Harnack inequality since it is a
key condition in our main results.

\begin{definition}
A function $h:M\rightarrow \mathbb{R}$ said to be \textit{harmonic%
\index{harmonic function}} on an open set $A\subset M$ \ if it is defined on
$%
\overline{A}$ and
\begin{equation}
h\left( x\right) =E_{x}\left( h\left( X_{T_{A}}\right) \right) \text{ \ for
all }x\in A.  \label{har1}
\end{equation}
\end{definition}

\begin{definition}
We will say that the \emph{elliptic Harnack inequality }$\left( \ref{H}%
\right) $\emph{\ }holds on $M$\emph{\ } if for all $x\in M,R>0$ and for any
non-negative harmonic function $u$ which is harmonic in $B(x,2R)$, the
following inequality holds
\begin{equation}
\sup\limits_{B\left( x,R\right) }u\leq H\inf\limits_{B\left( x,R\right) }u\,
\tag{H}  \label{H}
\end{equation}%
with some constant $H\geq 1$ independent of $x$ \ and $R.$
\end{definition}

The results of the present paper are the following.

\begin{theorem}
\label{tdistr} 1. If there is a $C>0$ such that the condition
$\left( \overline{E}\right) $
holds, that is there is a $C>0$ such that%
\begin{equation}
\sup_{y\in B\left( x,R\right) }E_{y}\left( x,R\right) \leq CE\left(
x,R\right)  \label{Ebar}
\end{equation}%
for all $x\in M,R>0$ then there is a \ $c>0$ \ such that for all $x\in
M,t,R>0$%
\begin{equation*}
\mathbb{P}_{x}\left( T_{x,R}<t\right) \leq \exp \left( -c\kappa \left(
x,t,R\right) \right)
\end{equation*}%
is true.\newline
\newline
2. If $M$ satisfies the elliptic Harnack inequality, then there are $b,C>0$
such that for all $x\in M,t,R>0$%
\begin{equation}
\mathbb{P}_{x}\left( T_{x,R}<t\right) \geq \exp \left( -C\nu \left(
x,t,bR\right) \right) .
\end{equation}
\end{theorem}

\begin{theorem}
\label{crrd}\label{tsta}Let us assume that there is an $R_{0}$ such that for
all $r<R_{0}$
\begin{equation*}
E\left( x,r\right) \simeq r^{\beta }
\end{equation*}%
holds with a $\beta >1$. Let $A,B\subset M$ be measurable sets $\ 0<\mu
\left( A\right) ,\mu \left( B\right) <\infty $. Then we have the upper part
of the short-time asymptotics:\
\begin{equation*}
\lim_{t\rightarrow 0}t^{\frac{1}{\beta -1}}\log P_{t}\left( A,B\right) \leq
-c\left[ d\left( A,B\right) \right] ^{\frac{\beta }{\beta -1}}
\end{equation*}%
\ and if in addition we assume that the $A,B$ sets are open and precompact
furthermore the elliptic Harnack inequality holds then
\begin{equation*}
\lim_{t\rightarrow 0}t^{\frac{1}{\beta -1}}\log P_{t}\left( A,B\right) \geq
-C\left[ d\left( A,B\right) \right] ^{\frac{\beta }{\beta -1}}.
\end{equation*}
\end{theorem}

\section{Discussion\protect\bigskip}

The usual sub-diffusive picture can be recovered assuming that $E\left(
x,R\right) \simeq R^{\beta }.$

\begin{corollary}
\label{c2}Let us assume that $\left( E_{\beta }\right) $ holds on $M$ for a $%
\beta >1$\newline
1. There are \ $c,C>0$ \ such that for all $x\in M,t,R>0,B=B\left(
x,R\right) $%
\begin{equation}
\mathbb{P}_{x}\left( T_{x,R}<t\right) \leq C\exp \left( -c\left[ \frac{%
R^{\beta }}{t}\right] ^{\frac{1}{\beta -1}}\right)  \label{uepsi}
\end{equation}%
is true.\newline
\newline
2. If $M$ satisfies the elliptic Harnack inequality, then there are $c,C>0$
such that
\begin{equation}
P\left( T_{x,R}<t\right) \geq c\exp \left( -C\left[ \frac{R^{\beta }}{t}%
\right] ^{\frac{1}{\beta -1}}\right) .  \label{lepsi}
\end{equation}
\end{corollary}

\begin{problem}
The classical bottle-neck construction shows that the condition $\left(
\overline{E}\right) $ does not imply the elliptic Harnack inequality. It
would be interesting to find an example which satisfies the elliptic Harnack
inequality but not $\left( \overline{E}\right) .$ At present we can not find
such one. \
\end{problem}

There are nice examples where the diffusion speed is \textquotedblright
direction\textquotedblright\ dependent (see \cite{B1} or \cite{HK}). We
briefly recall one following \cite{B1}. \ Consider the direct product $M=%
\mathbb{R}\times S_{2},$ where $S_{2}$ stands for the Sierpinski gasket and
let $Z_{t}=\left( X_{t},Y_{t}\right) $ be the process\ on it, where $X_{t}$
is the standard Wiener process, and $Y_{t}$ is the anomalous diffusion
process on $S_{2}$ (c.f. \cite{B2}) independent form $X_{t}.$ It is clear
that $X_{t}$ and $Y_{t}$ satisfy $\left( \ref{LEb}\right) ,\left( \ref{UEb}%
\right) $ with $\beta _{1}=2$ and $\beta _{2}>2$ respectively. Consequently
the diagonal upper estimate holds for both and for $Z_{t}$ \ as well while
the two-sided estimate is not true for any $\beta .$ It is then clear that
neither the elliptic Harnack inequality nor the short time asymtotics does
hold.

\section{Basic definitions}

\label{sdef}The Dirichlet form can be restricted to a set $A$ acting only on
functions with support \ in $A.$ \ The corresponding process is simply
killed on leaving $A$ (see \cite{FOT}). \ Let us denote the associated heat
kernel by $p_{t}^{A}\left( x,y\right) $ \ and the Green kernel by \ $%
g^{A}\left( x,y\right) .$

For sets we define
\begin{equation*}
d\left( A,B\right) =\inf_{x\in A,y\in B}d\left( x,y\right) .
\end{equation*}

To avoid technical difficulties we follow \cite{Ra1} and introduce
\begin{equation*}
P_{t}\left( A,B\right) =\int_{A}\int_{B}p_{t}\left( x,y\right) d\mu \left(
y\right) d\mu \left( x\right) .
\end{equation*}

\begin{definition}
We consider open metric balls defined by the metric $d(x,y)$ $x\in M,$ $R>0$
as
\begin{eqnarray*}
B(x,R) &=&\{y\in M:d(x,y)<R\}, \\
S(x,R) &=&\{y\in M:d(x,y)=R\}.
\end{eqnarray*}
\end{definition}

\begin{definition}
The exit time from a set $A$ is defined as
\begin{equation*}
T_{A}=\inf \{t>0:X_{t}\in A^{c}\},
\end{equation*}%
its expected value is denoted by
\begin{equation*}
E_{x}(A)=\mathbb{E}(T_{A}|X_{0}=x),
\end{equation*}%
and we will use the $E=E(x,R)=E_{x}(x,R)=E_{x}\left( B\left( x,R\right)
\right) $ short notations.
\end{definition}

\begin{definition}
The hitting time $\tau _{A}$ of a set $A$ is defined as the exit time of its
complement:
\begin{equation*}
\tau _{A}=T_{A^{c}}
\end{equation*}%
and for $A=B\left( x,R\right) $ we use the shorter form $\tau _{x,r}.$
\end{definition}

\begin{definition}
We introduce for a set $A\subset M,$
\begin{equation*}
\overline{E}\left( A\right) =\sup_{y\in A}E_{y}\left( A\right) .
\end{equation*}%
and for $x\in M,R>0$ we use the notation
\begin{equation*}
\overline{E}\left( x,R\right) =\overline{E}\left( B\left( x,R\right) \right)
.
\end{equation*}
\end{definition}

\begin{definition}
For any sets $A,B$ the capacity is defined via the Dirichlet form $\mathcal{E%
}$ \ \ by
\begin{equation*}
cap\left( A,B\right) =\inf \mathcal{E}\left( f,f\right) ,
\end{equation*}%
where the infimum runs for functions $f,$ $\ f|_{A}=1,f|_{B}=0.$ The
resistance is defined as
\begin{equation*}
\rho \left( A,B\right) =\frac{1}{cap\left( A,B\right) }.
\end{equation*}%
In particular we will use the following notations:
\begin{equation*}
\rho \left( x,r,R\right) =\rho \left( B\left( x,r\right) ,B^{c}\left(
x,R\right) \right) .
\end{equation*}
\end{definition}

\section{\protect\bigskip Distribution of the exit time}

In this section we show Theorem $\ref{tdistr}.$\bigskip\ First we recall a
result which was immediate from Lemma 5.3 of \cite{Td} for graphs and can be
seen in the same way for the present setup (see also \cite{B2}).\newline

\begin{proposition}
\label{p1}If we assume $\left( \overline{E}\right) $ then there is a $c>0$
such that for all $x\in M,r>0,t\leq \frac{1}{2}E\left( x,r\right) $%
\begin{equation}
\mathbb{P}_{x}\left( T_{x,r}>t\right) \geq c.  \label{ee}
\end{equation}
\end{proposition}

The proof of Theorem \ref{tdistr} based on the following observations. The probability of hitting a nearby ball in a "reasonable" time is bounded from
below if the elliptic Harnack inequality holds.

\begin{proposition}
\label{p2}If the elliptic Harnack inequality $\left( \ref{H}\right) $ holds then there are $c_{0},c_{1}>0$ such that
\begin{equation}
\mathbb{P}_{x}\left( \tau _{y,r}<s\right) \geq c_{0}.  \label{tt<}
\end{equation}%
provided $\frac{1}{4}d\left( x,y\right) \leq r\leq d\left( x,y\right) $ and $%
s>\frac{2}{c_{1}}E\left( x,9r\right) $
\end{proposition}

At this point we specify the constant $q$ which appears in the definition of $\kappa ,\nu $. \ Let $q=\frac{2}{c_{1}},$ which means, as we shall see in
the Lemma \ref{lptt>c}, that it depends via $c_{1}$ on the constant of the Harnack inequality.

The key observation is the following proposition. \ It provides a lower
bound for the probability hitting a ball in a given time.

\begin{proposition}
\label{p3}If the elliptic Harnack inequality $\left( \ref{H}\right) $ holds,
then there are $b,C>0$ such that for all $x,y\in M,t>0,\frac{1}{4}d\left(
x,y\right) \leq r<d=d\left( x,y\right) ,$%
\begin{equation}
\mathbb{P}_{x}\left( \tau _{y,r}<t\right) \geq \exp -C\nu \left(
x,t,bd\right) .  \label{cLE}
\end{equation}
\end{proposition}

\ First we need some lemmas.

\begin{lemma}
\label{lhg}If the elliptic Harnack inequality $\left( H\right) $ holds then
\ for $x\in M,R>r>0,B=B\left( x,R\right) ,A=B\left( x,r\right) $%
\begin{equation}
\inf_{w\in \overline{A}}g^{B}\left( w,x\right) \simeq \rho \left(
x,r,R\right) \simeq \sup_{w\in B\backslash A}g^{B}\left( w,x\right) .
\label{ehg}
\end{equation}
\end{lemma}

\begin{proof}
See Barlow's proof (\cite{Bnew}, Proposition 2) which generalizes
Propositions 4.1 and 4.3 of \cite{GT2} where the additional hypothesis of
bounded covering was used. Barlow's proof is given for weighted graphs, but
word by word the same proof works in the continuous case.
\end{proof}

\begin{lemma}
\label{lptt>c}If $M$ satisfies the elliptic Harnack inequality $\left( \ref%
{H}\right) $ then there is a $c_{1}>0$ such that for all $x\in M,r>0,w\in B=%
\overline{B}\left( x,4r\right) ,\frac{1}{4}d\left( w,x\right) \leq r$%
\begin{equation}
\mathbb{P}_{w}\left( \tau _{x,r}<T_{x,5r}\right) >c_{1}.  \label{tt>c}
\end{equation}
\end{lemma}

\begin{proof}
The investigated probability
\begin{equation}
u\left( w\right) =\mathbb{P}_{w}\left( \tau _{x,r}<T_{x,5r}\right)
\end{equation}%
is the capacity potential between $B^{c}\left( x,5r\right) $ and $B\left(
x,r\right) $ and clearly harmonic in $A=B\left( x,5r\right) \backslash
B\left( x,r\right) .$ \ So it can be decomposed in $A$%
\begin{equation*}
u\left( w\right) =\int_{\overline{B}\left( x,r\right) }g^{B\left(
x,5r\right) }\left( w,z\right) \pi \left( dz\right)
\end{equation*}%
where $\pi $ is the capacity measure with $\pi \left( \overline{B}\left(
x,r\right) \right) =1/\rho \left( x,r,5r\right) $ with support in $\overline{%
B}\left( x,r\right) $. \ From the maximum (minimum) principle it follows
that for any $z\in B\left( x,2r\right) $ the minimum of $g^{B\left(
x,5r\right) }\left( w,z\right) $ in $w\in \overline{B}\left( x,4r\right) $
is attained at a $w\in S\left( x,4r\right) $, (and the same applies for $%
u\left( w\right) $ as well). It follows then from the Harnack inequality for
$g^{B\left( x,5r\right) }\left( w,.\right) $ in $B\left( x,2r\right) $ that
for any $w\in S\left( x,4r\right) $
\begin{equation*}
\inf_{z\in \overline{B}\left( x,2r\right) }g^{B\left( x,5r\right) }\left(
w,z\right) \geq cg^{B\left( x,5r\right) }\left( w,x\right)
\end{equation*}%
\begin{equation*}
u\left( w\right) =\int_{\overline{B}\left( x,r\right) }g^{B\left(
x,5r\right) }\left( w,z\right) \pi \left( dz\right) \geq \frac{cg^{B\left(
x,5r\right) }\left( w,x\right) }{\rho \left( x,r,5r\right) }
\end{equation*}%
From Lemma \ref{lhg} we know that
\begin{equation*}
\sup_{y\in B\left( x,5r\right) \backslash B\left( x,4r\right) }g^{B\left(
x,5r\right) }\left( y,x\right) \simeq \rho \left( x,4r,5r\right) \simeq
\inf_{w\in \overline{B}\left( x,4r\right) }g^{B\left( x,5r\right) }\left(
w,x\right) .
\end{equation*}%
which means that
\begin{equation}
u\left( w\right) \geq c\frac{\rho \left( x,4r,5r\right) }{\rho \left(
x,r,5r\right) }  \label{ubig}
\end{equation}%
Similarly from Lemma \ref{lhg} it follows that
\begin{equation*}
\sup_{y\in B\left( x,5r\right) \backslash B\left( x,r\right) }g^{B\left(
x,5r\right) }\left( v,x\right) \simeq \rho \left( x,r,5r\right) \simeq
\inf_{w\in \overline{B}\left( x,r\right) }g^{B\left( x,5r\right) }\left(
w,x\right) .
\end{equation*}%
Finally if $y_{0}\in S\left( x,r\right) $ is on the ray from $x$ to $y\in
S\left( x,4r\right) $ then iterating the Harnack inequality along a finite
chain of balls of radius $r/4$ along this ray from $y_{0}$ to $y$ one obtains%
\begin{equation*}
g^{B\left( x,5r\right) }\left( y,x\right) \simeq g^{B\left( x,5r\right)
}\left( y_{0},x\right)
\end{equation*}%
which results that%
\begin{equation*}
\rho \left( x,r,5r\right) \leq c\rho \left( x,4r,5r\right) ,
\end{equation*}%
and the statement follows from $\left( \ref{ubig}\right) .$
\end{proof}

\begin{proof}[Proof of Proposition \protect\ref{p2}]
We insert the exit time $T_{x,9r}$ into the inequality $\tau _{y,r}<t$%
\begin{eqnarray*}
\mathbb{P}_{x}\left( \tau _{y,r}<t\right) &\geq &\mathbb{P}_{x}\left( \tau
_{y,r}<T_{x,9r}<t\right) \\
&=&\mathbb{P}_{x}\left( \tau _{y,r}<T_{x,9r}\right) -\mathbb{P}_{x}\left(
\tau _{y,r}<T_{x,9r},T_{x,9r}\geq t\right) \\
&\geq &\mathbb{P}_{x}\left( \tau _{y,r}<T_{x,9r}\right) -\mathbb{P}%
_{x}\left( T_{x,9r}\geq t\right) .
\end{eqnarray*}
On one hand the Markov inequality results that
\begin{equation*}
\mathbb{P}_{x}\left( T_{x,9r}\geq t\right) \leq \frac{E\left( x,9r\right) }{t%
}\leq \frac{E\left( x,9r\right) }{\frac{2}{c_{1}}E\left( x,9r\right) }%
<c_{1}/2
\end{equation*}
and on the other hand $B\left( y,5r\right) \subset B\left( x,9r\right) ,$
hence
\begin{equation*}
\mathbb{P}_{x}\left( \tau _{y,r}<T_{x,9r}\right) \geq \mathbb{P}_{x}\left(
\tau _{y,r}<T_{y,5r}\right) ,
\end{equation*}
and Lemma \ref{lptt>c} can be applied to get
\begin{equation*}
\mathbb{P}_{x}\left( \tau _{y,r}<T_{y,5r}\right) \geq c_{1}.
\end{equation*}
The result follows with $c_{0}=c_{1}/2.$
\end{proof}

\begin{lemma}
\label{lP>prod}For all $x\neq y\in M,t>0,l>1$
\begin{equation*}
\mathbb{P}_{x}\left( \tau _{y,r}<t\right) \geq \inf_{z,w\in \overline{B}%
\left( x,d\right) ,d\left( z,w\right) \leq 4r}\left[ \mathbb{P}_{z}\left(
\tau _{w,r}<s\right) \right] ^{l},
\end{equation*}%
where $d=d\left( x,y\right) ,$ $s=\frac{t}{l},r=\frac{d}{3l}$.
\end{lemma}

\begin{proof}
Let us consider a geodesic path $\pi $ from $x$ to $y$ Let $x_{1}\in \pi $
such that $d\left( x,x_{1}\right) =3r,$ $x_{2}\in \pi $ with $d\left(
x_{1},x_{2}\right) =3r$ etc. and finally $x_{l}=y.$ Let $\tau _{i}=\tau
_{x_{i},r}$ ,$\tau _{0}=0$ and $A_{i}=\left\{ \tau _{i}-\tau
_{i-1}<s\right\} $ for $i=1...l,\tau _{0}=0.$ Denote $A_{i}\left( z\right)
=\left\{ \tau _{i}-\tau _{i-1}<s,X_{\tau _{i}}=z\right\} $ and $\zeta
_{i}=X_{\tau _{i}}$. One can observe that $\prod_{i=0}^{l}A_{i}$ means that
the process spends less than time $s$ between the first hit of the
consecutive $B_{i}=B\left( x_{i},r\right) $ balls, consequently%
\begin{equation*}
\tau _{y,r}=\tau _{l}=\sum_{i=1}^{l}\tau _{i}-\tau _{i-1},
\end{equation*}%
\begin{equation*}
\left\{ \tau _{y,r}<t\right\} =\left\{ \sum_{i=1}^{l}\tau _{i}-\tau
_{i-1}<t\right\} \supset \prod_{i=0}^{l-1}A_{i}.
\end{equation*}%
Let us continue with the following estimates.
\begin{eqnarray*}
\mathbb{P}_{x}\left( \tau _{y,r}<t\right)  &\geq &\mathbb{P}_{x}\left(
\prod_{i=0}^{l}A_{i}\right) =\mathbb{E}_{x}\left( A_{l}\left( \zeta
_{l}\right) \prod_{i=0}^{l-1}A_{i}\right)  \\
&=&\mathbb{E}_{x}\left( I\left( A_{l}\right) |\prod_{i=0}^{l-1}I\left(
A_{i}\right) \right) \mathbb{P}_{x}\left( \prod_{i=0}^{l-1}I\left(
A_{i}\right) \right)
\end{eqnarray*}%
Now we use the strong Markov property (c.f \cite{FOT} $\left( A.2.3.\right)
^{\prime }$) to obtain the following.
\begin{eqnarray*}
&&\mathbb{E}_{x}\left( \mathbb{E}_{\zeta _{l-1}}\left[ I\left( A_{l}\right) %
\right] |\prod_{i=0}^{l-1}I\left( A_{i}\right) \right) \mathbb{P}_{x}\left(
\prod_{i=0}^{l-1}I\left( A_{i}\right) \right)  \\
&\geq &\inf_{z\in S\left( x_{l-1},r\right) }\mathbb{P}_{z}\left(
A_{l}\right) \mathbb{P}_{x}\left( \prod_{i=0}^{l-1}I\left( A_{i}\right)
\right)
\end{eqnarray*}%
Continuing the iteration one obtains%
\begin{eqnarray}
&&\mathbb{P}_{x}\left( \tau _{y,r}<t\right)   \notag \\
&\geq &\inf_{z\in S\left( x_{l-1},r\right) }\mathbb{P}_{z}\left(
A_{l}\right) \mathbb{P}_{x}\left( \prod_{i=1}^{l-1}I\left( A_{i}\right)
\right)   \notag \\
&\geq &\inf_{z,w\in B\left( x,d\right) ,d\left( z,w\right) \leq 4r}\left[
\mathbb{P}_{z}\left( \tau _{w,r}<s\right) \right] ^{l}.  \notag
\end{eqnarray}
\end{proof}

\begin{proof}[Proof of Proposition \protect\ref{p3}]
Let us apply Lemma \ref{lP>prod} with $l=\nu _{B}(t,3d)$\textbf{. \ }One
should observe that we have a uniform constant lower bound for $\mathbb{P}%
_{z}\left( \tau _{w,r}<s\right) $ by Proposition \ref{p2} provided
\begin{equation*}
s>\frac{2}{c_{1}}E\left( z,9r\right) .
\end{equation*}%
This condition is ensured by the definition of $\nu _{B}\left( t,3d\right) $
with $s=\frac{t}{\nu },r=\frac{d}{3l}$ since $9r=\frac{3d}{l}.$ Finally from
$l=\nu _{B\left( x,d\right) }(t,3d)\leq \nu _{B\left( x,3d\right)
}(t,3d)=\nu (x,t,3d)$ follows the statement.
\end{proof}

\begin{proof}[Proof of Theorem \protect\ref{tdistr}]
The upper estimate of Theorem \ref{tdistr} can be seen along the
lines of the proof of Theorem 5.1 in \cite{Td}. The lower bound is immediate
from $\left( \ref{cLE}\right) $. Let $b=6$ then for any $y\in S\left(
x,2R\right) ,r=\frac{6R}{l},l=\nu \left( x,t,6R\right) $%
\begin{equation*}
\mathbb{P}_{x}\left( T_{x,R}<t\right) \geq \mathbb{P}_{x}\left( \tau
_{y,r}<t\right)
\end{equation*}%
and the result follows from Proposition \ref{p3}.
\end{proof}

\section{Short time asymptotics}

Before we start the proof let us make it clear why we assume that $E\left(
x,r\right) \simeq r^{\beta }$ for $r<R_{0}$. \ Without this assumption the
upper and lower estimates will not meet. \ Using the proof which will
follow, one can obtain separate upper and lower bound for the short time
asymptotics. \ Let us mention that as $t$ goes to zero the functions $\kappa
$ and $\nu $ go to infinity (c.f. Definition \ref{dKL}) consequently for
enough small $t$ $\frac{r}{\kappa }$ and $\frac{r}{\nu }$ fall below $R_{0}$
and the polynomial approximation can be used and provides us the more
readable and coinciding upper and lower bound.

\begin{proof}[Proof of Theorem \protect\ref{tsta}]
\ Consider $A,B\subset M$, denote $d=d\left( A,B\right) $ and let us use $%
\left( \ref{uepsi}\right) $ to get
\begin{eqnarray*}
P_{t}\left( A,B\right)  &=&\int_{A}\left( P_{t}1_{B}\right) \left( x\right)
d\mu \left( x\right)  \\
&\leq &C\mu \left( A\right) \mu \left( B\right) \exp \left[ -c\left( \frac{%
d^{\beta }}{t}\right) ^{\frac{1}{\beta -1}}\right] ,
\end{eqnarray*}%
which results that
\begin{equation*}
t^{\frac{1}{\beta -1}}\log P_{t}\left( A,B\right) \leq t^{\frac{1}{\beta -1}}%
\left[ C+\log \mu \left( A\right) +\log \mu \left( B\right) \right]
-cd\left( A,B\right) ^{\frac{\beta }{\beta -1}},
\end{equation*}%
and
\begin{equation*}
\lim_{t\rightarrow 0}t^{\frac{1}{\beta -1}}\log \left( P_{t}\left(
A,B\right) \right) \leq -cd\left( A,B\right) ^{\frac{\beta }{\beta -1}}.
\end{equation*}%
\ If $A,B$ are open we can find for any $x\in A,y\in B$ a very small $r$
such that the balls $A_{r}=B\left( x,r\right) \subset A,$ $B_{r}=B\left(
y,r\right) \subset B.$ For the lower estimate we decompose the path with the
first hit of $B\left( y,r/2\right) \subset B\left( y,r\right) \subset B.$
The only task is to show that the probability that the process stays in $%
B\left( y,r\right) $ until $t$ is bounded from below by a constant. \ Let us
observe that
\begin{equation*}
R=d\left( x,y\right) \geq d\left( A,B\right) +2r.
\end{equation*}%
Denote $\xi =X_{\tau },$where $\tau =\tau _{y,r/2}$,
\begin{equation*}
P_{t}\left( A,B\right) =\int_{A}P_{t}\left( x,B\right) d\mu \left( x\right)
\geq \int_{A}P_{t}\left( x,B_{r}\right) d\mu \left( x\right) .
\end{equation*}%
For any fixed $x\in A$
\begin{equation*}
P_{t}\left( x,B_{r}\right) \geq \mathbb{E}_{x}\left( I\left( \tau <t\right)
I\left( X_{t}\in B_{r}\right) \right)
\end{equation*}%
Now we use the strong Markov property(c.f \cite{FOT} $\left( A.2.3.\right)
^{\prime }$) to decompose the path according to $\xi =X_{\tau }$ and obtain
\begin{eqnarray*}
&&\mathbb{E}_{x}\left( I\left( \tau <t\right) I\left( X_{t}\in B_{r}\right)
\right)  \\
&=&\mathbb{E}_{x}\left( I\left( \tau <t\right) \mathbb{E}_{\xi }\left[
I\left( X_{t-\tau }\in B_{r}\right) \right] \right)  \\
&\geq &\mathbb{E}_{x}\left( I\left( \tau <t\right) \mathbb{E}_{\xi }\left[
I\left( T_{y,r}>t-\tau \right) \right] \right)  \\
&\geq &\mathbb{E}_{x}\left( I\left( \tau <t\right) \mathbb{E}_{\xi }\left[
I\left( T_{\xi ,r/2}>t\right) \right] \right)  \\
&\geq &\mathbb{E}_{x}\left( I\left( \tau <t\right) \inf\limits_{w\in S\left(
y,r/2\right) }\mathbb{E}_{w}\left[ I\left( T_{w,r/2}>t\right) \right]
\right) ,
\end{eqnarray*}%
It is clear from $\left( \ref{ee}\right) $ that for all $w$
\begin{equation*}
\mathbb{P}_{w}\left( T_{w,r/2}>t\right) \geq c^{\prime }
\end{equation*}%
if
\begin{equation*}
t<\frac{1}{2}E\left( w,r/2\right) .
\end{equation*}%
This can be ensured using the lower bound of $\left( E_{\beta }\right) $ for
small $t:$%
\begin{equation}
t<\frac{c}{2}\left( \frac{r}{2}\right) ^{\beta }<\frac{1}{2}E\left(
w,r/2\right) .  \label{cond1}
\end{equation}%
So if $\left( \ref{cond1}\right) $ holds we have
\begin{eqnarray*}
P_{t}\left( x,B\left( y,r\right) \right)  &\geq &c^{\prime }\mathbb{E}%
_{x}\left( I\left( \tau <t\right) \right)  \\
&=&c^{\prime }\mathbb{P}_{x}\left( \tau <t\right) .
\end{eqnarray*}%
Set $\rho =\varepsilon r,B^{\prime }=B\left( y,\rho \right) $. The next step
is to use $\left( \ref{cLE}\right) $ {\small \ to get}
\begin{equation*}
\mathbb{P}_{x}\left( \tau <t\right) \geq \exp \left( -C\left( \frac{R^{\beta
}}{t}\right) ^{\frac{1}{\beta -1}}\right) .
\end{equation*}%
The proper choice of the constants follows from the restrictions:
\begin{equation}
\varepsilon r=\rho ,\varepsilon <1/2,  \label{e1}
\end{equation}%
\
\begin{equation}
\frac{t}{k}\simeq \rho ^{\beta }=\left( \frac{R}{3k}\right) ^{\beta },
\label{e2}
\end{equation}%
and
\begin{equation}
t^{\frac{1}{\beta }}<cr.  \label{e3}
\end{equation}%
If we consider
\begin{equation*}
\rho =\frac{R}{3k}=\varepsilon r
\end{equation*}%
it follows that the proper choice for $\varepsilon $ \ is
\begin{equation*}
\varepsilon =c\left( \frac{r}{R}\right) ^{\frac{1}{\beta -1}}.
\end{equation*}%
The short time asymptotics now is immediate.
\begin{equation*}
P_{t}\left( A,B\right) \geq \mu \left( A\right) c\exp \left( -C\left( \frac{%
R^{\beta }}{t}\right) ^{\frac{1}{\beta -1}}\right) ,
\end{equation*}%
\begin{equation*}
\lim\limits_{t\rightarrow 0}t^{\frac{1}{\beta -1}}\log P_{t}\left(
A,B\right) \geq -Cd^{\frac{\beta }{\beta -1}}\left( x,y\right) .
\end{equation*}%
Finally we let $d\left( x,y\right) \rightarrow d\left( A,B\right) $ and we
receive the lower bound.
\end{proof}

\end{document}